# Bat Algorithm: A Novel Approach for Global Engineering Optimization

Xin-She Yang

Mathematics and Scientific Computing, National Physical Lab, Teddington TW11 0LW, UK

Email: xy227@cam.ac.uk

Amir Hossein Gandomi

Department of Civil Engineering, The University of Akron, Akron, OH 44325, USA

Email: a.h.gandomi@gmail.com

**Abstract**

**Purpose** – Nature-inspired algorithms are among the most powerful algorithms for optimization. In this study, a new nature-inspired metaheuristic optimization algorithm, called bat algorithm (BA), is introduced for solving engineering optimization tasks.

**Design/methodology/approach** – The proposed BA is based on the echolocation behavior of bats. After a detailed formulation and explanation of its implementation, BA is verified using eight nonlinear engineering optimization problems reported in the specialized literature.

**Findings** – BA has been carefully implemented and carried out optimization for eight well-known optimization tasks. Then, a comparison has been made between the proposed algorithm and other existing algorithms.

**Practical implications** – Engineering Design Optimization

**Originality/value** – The optimal solutions obtained by the proposed algorithm are better than the best solutions obtained by the existing methods. The unique search features used in BA are analyzed, and their implications for future research are also discussed in detail.

# 1. Introduction

Design optimization forms an important part of any design problem in engineering and industry. Structural design optimization focuses on finding the optimal and practical solutions to complex structural design problems under dynamic complex loading pattern with complex nonlinear constraints. These constraints often involve thousands of and even millions of members with stringent limitations on stress, geometry as well as loading and service requirements. The aim is not only to minimize the cost and materials usage, but also to maximize their performance and lifetime service. All these designs are of scientific and practical importance (Deb 1995, Yang 2010). However, most structural design optimization problems are highly nonlinear and multimodal with noise, and thus they are often NP-hard. Finding the right and practically efficient algorithms are usually difficult, if not impossible. In realistic, the choice of an algorithm requires extensive experience and knowledge of the problem of interest. Even so, there is no guarantee that an optimal or even suboptimal solution can be found.

Metaheuristic algorithms including evolutionary and swarm intelligence algorithms are now becoming powerful methods for solving many tough problems (Gandomi and Alavi 2011) and especially real-world engineering problems (Gandomi et al. 2011, Alavi and Gandomi 2011). The vast majority of heuristic and metaheuristic algorithms have been derived from the behavior of biological systems and/or physical systems in nature. For example, particle swarm optimization was developed based on the swarm behavior of birds and fish (Kennedy and Eberhart 1995) or charged system search inspired from physical processes (Kaveh and Talatahari 2010). New algorithms are also emerging recently, including harmony search and the firefly algorithm. The former was inspired by the improvising process of composing a piece of music (Geem et al. 2001), while the latter was formulated based on the flashing behavior of fireflies (Yang 2008). Each of these algorithms has certain advantages and disadvantages. For example, simulating annealing (Kirkpatrick et al. 1983) can almost guarantee to find the optimal solution if the cooling process is slow enough and the simulation is running long enough; however, the fine adjustment in parameters does affect the convergence rate of the optimization process. A natural question is whether it is possible to combine major advantages of these algorithms and try to develop a potentially better algorithm?

This paper is such an attempt to address this issue. In this paper, we intend to propose a new metaheuristic method, namely, the Bat Algorithm (BA), based on the echolocation behavior of bats, and preliminary studies show that this algorithm is very promising (Yang 2010). The capability of echolocation of microbats is fascinating as these bats can find their prey and discriminate different types of insects even in complete darkness. We will first formulate the bat algorithm by idealizing the echolocation behavior of bats. We then describe how it works and make comparison with other existing algorithms. Finally, we will discuss some implications for further studies.

**2 Echolocation of Microbats**

Bats are fascinating animals. They are the only mammals with wings and they also have advanced capability of echolocation. It is estimated that there are about 1000 different species which account for up to about one fifth of all mammal species (Altringham 1996). Their size ranges from the tiny bumblebee bat (of about 1.5 to 2 g) to the giant bats with wingspan of about 2 m and weight up to about 1 kg. Microbats typically have forearm length of about 2.2 to 11 cm. Most bats uses echolocation to a certain degree; among all the species, microbats are a famous example as microbats use echolocation extensively while megabats do not (Richardson 2008).

Most microbats are insectivores. Microbats use a type of sonar, called, echolocation, to detect prey, avoid obstacles, and locate their roosting crevices in the dark. These bats emit a very loud sound pulse and listen for the echo that bounces back from the surrounding objects. Their pulses vary in properties and can be correlated with their hunting strategies, depending on the species. Most bats use short, frequency-modulated signals to sweep through about an octave, while others more often use constant-frequency signals for echolocation. Their signal bandwidth varies depends on the species, and often increased by using more harmonics.

Though each pulse only lasts a few thousandths of a second (up to about 8 to 10 ms); however, it has a constant frequency which is usually in the region of 25 kHz to 150 kHz. The typical range of frequencies for most bat species are in the region between 25 kHz and 100 kHz, though some species can emit higher frequencies up to 150 kHz. Each ultrasonic burst may last typically 5 to 20 ms, and microbats emit about 10 to 20 such sound bursts every second. When hunting for prey, the rate of pulse emission can be sped up to about 200 pulses per second when they fly near their prey. Such short sound bursts imply the fantastic ability of the signal processing power of bats. In fact, studies show the integration time of the bat ear is typically about 300 to 400 μs.

As the speed of sound in air is typically $v$ = 340 m/s, the wavelength $\lambda$ of the ultrasonic sound bursts with a constant frequency f is given by $\lambda = v/f$, which is in the range of 2 mm to 14 mm for the typical frequency range from 25 kHz to 150 kHz. Such wavelengths are in the same order of their prey sizes.

Amazingly, the emitted pulse could be as loud as 110 dB, and, fortunately, they are in the ultrasonic region. The loudness also varies from the loudest when searching for prey and to a quieter base when homing towards the prey. The travelling range of such short pulses is typically a few meters, depending on the actual frequencies (Richardson 2008). Microbats can manage to avoid obstacles as small as thin human hairs.

Studies show that microbats use the time delay from the emission and detection of the echo, the time difference between their two ears, and the loudness variations of the echoes to build up three dimensional scenario of the surrounding. They can detect the distance and orientation of the target, the type of prey, and even the moving speed of the prey such as small insects. Indeed, studies suggested that bats seem to be able to discriminate targets by the variations of the Doppler Effect induced by the wing-flutter rates of the target insects (Altringham 1996).

Obviously, some bats have good eyesight, and most bats also have very sensitive smell sense. In reality, they will use all the senses as a combination to maximize the efficient detection of prey and smooth navigation. However, here we are only interested in the echolocation and the associated behavior.

Such echolocation behavior of microbats can be formulated in such a way that it can be associated with the objective function to be optimized, and this makes it possible to formulate new optimization algorithms. In the rest of this paper, we will first outline the basic formulation of the Bat Algorithm (BA) and then discuss the implementation and comparison in detail.

**3 Bat Algorithm**

If we idealize some of the echolocation characteristics of microbats, we can develop various bat-inspired algorithms or bat algorithms. For simplicity, we now use the following approximate or idealized rules:

1. All bats use echolocation to sense distance, and they also `know' the difference between food/prey and background barriers in some magical way;
2. Bats fly randomly with velocity $v_i$ at position $x_i$ with a fixed frequency $f_{min}$, varying wavelength $\lambda$ and loudness $A_0$ to search for prey. They can automatically adjust the wavelength (or frequency) of their emitted pulses and adjust the rate of pulse emission $r$ in the range of [0, 1], depending on the proximity of their target;

3. Although the loudness can vary in many ways, we assume that the loudness varies from a large (positive) $A_0$ to a minimum constant value $A_{min}$.

Another obvious simplification is that no ray tracing is used in estimating the time delay and three dimensional topography. Though this might be a good feature for the application in computational geometry, however, we will not use this feature, as it is more computationally extensive in multidimensional cases.

In addition to these simplified assumptions, we also use the following approximations, for simplicity. In general the frequency f in a range [$f_{min}$, $f_{max}$] corresponds to a range of wavelengths [$\lambda_{min}$, $\lambda_{max}$]. For example a frequency range of [20 kHz, 500 kHz] corresponds to a range of wavelengths from 0.7 mm to 17 mm.

For a given problem, we can also use any wavelength for the ease of implementation. In the actual implementation, we can adjust the range by adjusting the wavelengths (or frequencies), and the detectable range (or the largest wavelength) should be chosen such that it is comparable to the size of the domain of interest, and then toning down to smaller ranges. Furthermore, we do not necessarily have to use the wavelengths themselves; instead, we can also vary the frequency while fixing the wavelength $\lambda$. This is because $\lambda$ and $f$ are related due to the fact $\lambda f$ is constant. We will use this later approach in our implementation.

For simplicity, we can assume f is within [0, $f_{max}$]. We know that higher frequencies have short wavelengths and travel a shorter distance. For bats, the typical ranges are a few meters. The rate of pulse can simply be in the range of [0, 1] where 0 means no pulses at all, and 1 means the maximum rate of pulse emission.

Based on these approximations and idealization, the basic steps of the Bat Algorithm (BA) can be summarized as the pseudo code shown in Fig. 1.

**Fig. 1** Pseudo code of the bat algorithm (BA)

### 3.1 Velocity and Position Vectors of Virtual Bats

In simulations, we use virtual bats naturally. We have to define the rules how their positions $x_i$ and velocities $v_i$ in a d-dimensional search space are updated. The new solutions $x^t_i$ and velocities $v^t_i$ at time step $t$ are given by

$$f_i = f_{min} + (f_{max} - f_{min})\beta \quad (1)$$
$$v^t_i = v^{t-1}_i + (x^t_i - x_*)f_i \quad (2)$$
$$x^t_i = x^{t-1}_i + v^t_i \quad (3)$$

where $\beta \in [0,1]$ is a random vector drawn from a uniform distribution. Here $x_*$ is the current globalbest location (solution) which is located after comparing all the solutions among all the n bats. As the product $\lambda_i f_i$ is the velocity increment, we can use either $f_i$ (or $\lambda_i$) to adjust the velocity change while fixing the other factor $\lambda_i$ (or $f_i$), depending on the type of the problem of interest. In our implementation, we will use $f_{min} = 0$ and $f_{max} = 100$, depending the domain size of the problem of interest. Initially, each bat is randomly assigned a frequency that is drawn uniformly from $[f_{min}, f_{max}]$.

For the local search part, once a solution is selected among the current best solutions, a new solution for each bat is generated locally using a local random walk:

$$x_{new} = x_{old} + \varepsilon A^t \tag{4}$$

where $\varepsilon \in [-1,1]$ is a random number, while $A^t = <A^t_i>$ is the average loudness of all the bats at this time step.

The update of the velocities and positions of bats have some similarity to the procedure in the standard particle swarm optimization (Geem et al. 2001) as $f_i$ essentially controls the pace and range of the movement of the swarming particles. To a degree, BA can be considered as a balanced combination of the standard particle swarm optimization and the intensive local search controlled by the loudness and pulse rate.

### 3.2 Variations of Loudness and Pulse Emission

Furthermore, the loudness $A_i$ and the rate $r_i$ of pulse emission have to be updated accordingly as the iterations proceed. As the loudness usually decreases once a bat has found its prey, while the rate of pulse emission increases, the loudness can be chosen as any value of convenience. For example, we can use $A_0 = 100$ and $A_{min} = 1$. For simplicity, we can also use $A_0 = 1$ and $A_{min} = 0$, assuming $A_{min} = 0$ means that a bat has just found the prey and temporarily stop emitting any sound. Now we have

$$A^{t+1}_i = \alpha A^t_i, \quad r^{t+1}_i = r_i^0 \ [1- \exp(-\gamma t)] \tag{5}$$

where $\alpha$ and $\gamma$ are constants. In fact, $\alpha$ is similar to the cooling factor of a cooling schedule in the simulated annealing (Kirkpatrick et al. 1983, Yang 2008). For any $0 < \alpha < 1, 0 < \gamma$, we have

$$A^t_i \to 0, \quad r^t_i \to r_i^0, \qquad as\ t \to \infty \tag{6}$$

In the simplicity case, we can use $\alpha = \gamma$, and we have in fact used $\alpha = \gamma = 0.9$ in our simulations. The choice of parameters requires some experimenting. Initially, each bat should have different values of loudness and pulse emission rate, and this can be achieved by randomization. For example, the initial loudness $A^0_i$ can typically be [1, 2], while the initial emission rate $r^0_i$ can be around zero, or any value $r^0_i \in [0, 1]$ if using Eq. 5. Their loudness and emission rates will be updated only if the new solutions are improved, which means that these bats are moving towards the optimal solution.

## 4. Non-linear Engineering Design Tasks

Most real-world engineering optimization problems are nonlinear with complex constraints, sometimes the optimal solutions of interest do not even exist. In order to see how BA performs, we will now test it against some well-known, tough but yet diverse, benchmark design problems. We have chosen eight case studies as:

- mathematical problem,
- Himmelblau's problem,
- three-bar truss design,
- speed reducer design,
- parameter identification of structures,
- cantilever stepped beam,
- heater exchanger design, and
- car side problem

The reason for such choice is to provide a validation and test of the proposed BA against a diverse range of real-world engineering optimization problems. As we will see below, for most problems, the optimal solutions obtained by BA are far better than the best solutions reported in the literature. In all case studies, the statistical measures have been obtained, based on 50 independent runs.

### 4.1. Case 1: Mathematical Problem

Now let us start with a nonlinear mathematical benchmark problem. This problem has been used as a benchmark constrained optimization problem with some active inequality constraints (Chen and Vassiliadis 2003). In this problem, $N$ is the number of variables and it is a multiple of four ($N=4n$, $n=1, 2, 3, ...$). This problem has $N/2$ inequality constraints and $2N$ simple bounds or limits. The problem can be stated as follows:

Minimize: $f(X) = \sum_{i=1}^{N} \sqrt{i}(x_i - 1)^2 + \left[\sum_{i=1}^{N} x_i^2 - 25\right]^2$ (7)

Subject to:

$g_j = x_{4(j-1)+1} + 2x_{4(j-1)+2} + 3x_{4(j-1)+3} + 4x_{4(j-1)+4} - 20$ (8)

$$0 \leq g_j \leq 30 \quad \text{and} \quad j = 1, 2, \ldots, \frac{N}{4} \tag{9}$$

with simple bounds

$0.5 \leq x_i \leq 10$ ($i=1,2,\ldots N$).

For this problem, the global optimum and best known optimum for *N*=12 and *N*=60 obtained by BA is given in Table 1. It can clearly be seen from Table 1 that the BA successfully find the global minimum. The statistical results of the mathematical problem are also presented in Table 2.

**Table 1** comparison of the BA results and global optimums for the mathematical problem

**Table 2** statistical results of the mathematical problem

### 4.2. Case 2: Himmelblau's Problem

Now we solve a well-known benchmark problem, namely Himmelblau's problem. This problem was originally proposed by Himmelblau's (Himmelblau 1972) and it has been widely used as a benchmark nonlinear constrained optimization problem. In this problem, there are five design variables [$x_1, x_2, x_3, x_4, x_5$], six nonlinear inequality constraints, and ten simple bounds or limits. The problem can be stated as follows:

Minimize: $f(X) = 5.3578547 x_3^2 + 0.8356891 x_1 x_5 + 37.293239 x_1 - 40792.141$ (10)

Subject to $0 \leq g_1 \leq 92$, $90 \leq g_2 \leq 110$, and $20 \leq g_3 \leq 25$ where

$$g_1 = 85.334407 + 0.0056858 x_1 x_5 + 0.0006262 x_1 x_4 - 0.0022053 x_3 x_5 \tag{11}$$

$$g_2 = 80.51249 + 0.0071317 x_2 x_5 + 0.0029955 x_1 x_2 - 0.0021813 x_3^2 \tag{12}$$

$$g_3 = 9.300961 + 0.0047026 x_3 x_5 + 0.0012547 x_1 x_3 - 0.0019085 x_3 x_4 \tag{13}$$

with simple bounds

$78 \leq x_1 \leq 102$, $33 \leq x_2 \leq 45$, and $27 \leq x_3, x_4, x_5 \leq 45$.

The best known optimum for the Himmelblau's problem obtained by BA is given in Table 3.

**Table 3** BA results for the Himmelblau's problem.

The problem was initially solved by Himmelblau (Himmelblau 1972) using a generalized gradient method. Since then, this problem has also been solved using several other methods such as GA (Gen and Cheng 1997, Homaifar et al. 1994), harmony search (HS) algorithm (Lee and Geem 2004, Fesanghary et al. 2008), and PSO (He et al. 2004, Shi and Eberhart 1998). Table II summarizes the results obtained by BA, as well as those published in the literature. It can clearly be seen from Table 4 that the result obtained by BA is better than the best feasible solution previously reported.

**Table 4** Statistical results for the Himmelblau's problem

### 4.3. Case 3: A Three Bar Truss Design

This case study considers a 3-bar planar truss structure shown in Fig. 2. This problem was first presented by Nowcki (1974). The volume of a statically loaded 3-bar truss is to be minimized subject to stress (σ) constraints on each of the truss members. The objective is to evaluate the optimal cross sectional areas. The mathematical formulation is given as below:

$$\text{Minimize: } f(X) = \left(2\sqrt{2x_1} + x_2\right) \times l \quad (14)$$

Subject to:

$$g_1 = \frac{\sqrt{2}x_1 + x_2}{\sqrt{2}x_1^2 + 2x_1 x_2} P - \sigma \leq 0 \quad (15)$$

$$g_2 = \frac{x_2}{\sqrt{2}x_1^2 + 2x_1 x_2} P - \sigma \leq 0 \quad (16)$$

$$g_3 = \frac{1}{x_1 + \sqrt{2}x_2} P - \sigma \leq 0 \quad (17)$$

where

$$0 \leq x_1 \leq 1 \text{ and } 0 \leq x_2 \leq 1; l = 100\,cm, P = 2KN/cm^2, \text{and } \sigma = 2KN/cm^2$$

**Fig. 2** Three-bar truss.

This design problem is a nonlinear fractional programming problem. The statistical values of the best solution obtained by BA are given in Table 5. The best solution by BA is ($x_1$, $x_2$) = (0.78863, 0.40838) with the objective value equal to 263.896248. Table 6 presents the best solutions obtained by BA and those reported by Ray and Saini (2001) and Tsai (2005). It can be seen clearly that the best objective value reported by Tsai (2005) is not feasible because the first constraint ($g_1$) is violated. Hence, it can be concluded that the results obtained by BA are better than those of the previous studies.

**Table 5** Statistical results of the best three bar truss model

**Table 6** Best solutions for the three bar truss design example

### 4.4. Case 4: Speed Reducer Design

The design of a speed reducer is a more complex case study (Golinski 1973) and it is one of the benchmark structural engineering problems (Gandomi and Yang 2011). This problem involves seven design variables, as shown in Fig. 3, with the face width $b$ ($x_1$), module of teeth $m$ ($x_2$), number of teeth on pinion $z$ ($x_3$), length of first shaft between bearings $l_1$ ($x_4$), length of second shaft between bearings $l_2$ ($x_5$), diameter of first shaft $d_1$ ($x_6$), and diameter of second shaft $d_2$ ($x_7$). The objective is to minimize the total weight of the speed reducer. There are nine constraints, including the limits on the bending stress of the gear teeth, surface stress, transverse deflections of shafts 1 and 2 due to transmitted force, and stresses in shafts 1 and 2.

**Fig. 3** Speed reducer.

The mathematical formulation can be summarized as follows:

Minimize:
$$f(X) = 0.7854 x_1 x_2^2 (3.3333 x_3^2 + 14.9334 x_3 - 43.0934) - 1.508 x_1 (x_6^2 + x_7^2) \\ + 7.477 (x_6^3 + x_7^3) + 0.7854 (x_4 x_6^2 + x_5 x_7^2)$$

(21)

Subject to:

$$g_1 = \frac{27}{x_1 x_2^2 x_3} P - 1 \leq 0 \tag{18}$$

$$g_2 = \frac{397.5}{x_1 x_2^2 x_3^2} - 1 \leq 0 \tag{19}$$

$$g_3 = \frac{1.93}{x_2 x_3 x_4^3 x_6^4} - 1 \leq 0 \tag{20}$$

$$g_4 = \frac{1.93}{x_2 x_3 x_4^3 x_7^4} - 1 \leq 0 \tag{21}$$

$$g_5 = \frac{\sqrt{\left(\frac{745 x_4}{x_2 x_3}\right)^2 + 1.69 \times 10^6}}{110 x_6^3} - 1 \leq 0 \tag{22}$$

$$g_6 = \frac{\sqrt{\left(\frac{745 x_4}{x_2 x_3}\right)^2 + 157.5 \times 10^6}}{85 x_7^3} - 1 \leq 0 \tag{23}$$

$$g_7 = \frac{x_2 x_3}{40} - 1 \leq 0 \tag{24}$$

$$g_8 = \frac{5 x_2}{B - 1} - 1 \leq 0 \tag{25}$$

$$g_9 = \frac{x_1}{12 x_2} - 1 \leq 0 \tag{26}$$

In addition, the design variables are also subject to simple bounds list in Table VII. This problem has been solved by using BA, and the corresponding statistical values of the best solutions are also presented in Table 7.

**Table 7** Statistical results of the speed reducer design example

Table 8 summarizes a comparison of the results obtained by BA with those obtained by other methods. Although some of the best objective values are better than those of BA, these reported values are not feasible because some of the constraints are violated. Thus, BA obtained the best feasible solution for this problem.

**Table 8** Statistical results of the speed reducer design example

### 4.5. Case 5: Parameter Identification of Structures

Estimation of structural parameter is the art of reconciling an *a priori* finite-element model (FEM) of the structure with nondestructive test data. It has a great potential for use in FEM updating. Saltenik and Sanayei (1996) developed a parameter estimation benchmark using measured strains for simultaneous estimation of the structural parameters. The parameter estimation objective function is defined as follows:

$$\text{Minimize:} \sum_{i=1}^{NMS} \left| \frac{\left([\varepsilon_a]_{m,i} - [\varepsilon_a]_{a,i}\right)}{[\varepsilon]_{m,i}} \right| \tag{27}$$

where $[\varepsilon_a]_m$ is the measured strains, $[\varepsilon_a]_m$ = number of measurements (*NMS*) × number of loading states (*NLS*), and $[\varepsilon_a]_a$ is the analytical strains.

The static FEM equation for a structural system is $[F]=[K][U]$. Thus, the analytical strains can be calculated as follows:

$$[\varepsilon] = [B][K]^{-1}[F] \tag{28}$$

It is not required to measure all the strains, therefore, Eq. (32) is partitioned based on measured strain ***a*** and unmeasured strain ***b***:

$$\begin{bmatrix} \varepsilon_a \\ \varepsilon_b \end{bmatrix} = \begin{bmatrix} B_a \\ B_b \end{bmatrix} [K]^{-1}[F] \tag{29}$$

Since there is no need for unmeasured strains $[\varepsilon_b]$ is eliminated as:

$$[\varepsilon_a] = [B_a][K]^{-1}[F] \tag{30}$$

In this work, the case study is a frame structure presented by Saltenik and Sanayei (1996) (see Fig. 4). The identified parameter in this example is moment of inertia *I* (***X***) for each member.

**Fig. 4** Frame structure used for parameter identification example.

A 445 N load is applied to degrees of freedom of 2, 5, 8 and 11, and each load set is composed of only one force. Strains are measured on 3, 6 and 7 for each load set. The cross section areas are respectively 484 cm$^2$ and 968 cm$^2$ for the horizontal and inclined members. The Elastic modulus is 206.8 GPa for all elements. The optimal solution is obtained at ***X*** = [869, 869, 869, 869, 869, 1320, 1320] (cm$^4$) with corresponding function value equal to *f\**(***X***) = 0.00000. The statistical results for this case study provided by BA are presented in Table 9.

> **Table 9** Best solutions for the parameter identification example using BA

The analytical algorithm proposed by Saltenik and Sanayei (1996) is not applicable to this problem due to a singularity. Arjmandi (2010) solved this problem using GA. A comparison of the results obtained by GA and BA with the measured values is summarized in Fig. 5. The results show that BA has found the global optimum and identified all the parameters without any error.

> **Fig. 5** Parameter identification results using GA and BA

### 4.6. Case 6: Cantilever Stepped Beam

The capability of BA for continuous and discrete variable design problems are verified using a design problem with ten variables. The case is originally presented by Thander and Vanderplaates (1995). Fig. 6 presents a five-stepped cantilever beam with rectangular shape. In this case study, the width ($x_1$-$x_5$) and height ($x_6$-$x_{10}$) of the beam in all five steps of the cantilever beam are design variables. The volume of the beam is to be minimized. The objective function is formulated as follows:

Minimize:  $$V = \sum_{i=1}^{5} x_i x_{i+5} l_i \qquad (31)$$

where $l_i$ = 100 cm (i=1,2,..,5)

> **Fig. 6** A stepped cantilever beam.

Subject to the following constraints:

$$g_1 = \frac{600P}{x_5 x_{10}^2} - 14000 \leq 0 \qquad (32)$$

$$g_2 = \frac{6P(l_s + l_4)}{x_4 x_9^2} - 14000 \leq 0 \qquad (33)$$

$$g_3 = \frac{6P(l_s + l_4 + l_3)}{x_3 x_8^2} - 14000 \leq 0 \qquad (34)$$

$$g_4 = \frac{6P(l_s + l_4 + l_3 + l_2)}{x_2 x_7^2} - 14000 \leq 0 \tag{35}$$

$$g_5 = \frac{6P(l_s + l_4 + l_3 + l_2 + l_1)}{x_1 x_6^2} - 14000 \leq 0 \tag{36}$$

$$g_6 = \frac{Pl^3}{3E}\left(\frac{1}{I_s} + \frac{7}{I_4} + \frac{19}{I_3} + \frac{37}{I_2} + \frac{61}{I_1}\right) - 2.7 \leq 0 \tag{37}$$

$$g_7 = \frac{x_{10}}{x_5} - 20 \leq 0 \tag{38}$$

$$g_8 = \frac{x_9}{x_4} - 20 \leq 0 \tag{39}$$

$$g_9 = \frac{x_8}{x_3} - 20 \leq 0 \tag{40}$$

$$g_{10} = \frac{x_7}{x_2} - 20 \leq 0 \tag{41}$$

$$g_{11} = \frac{x_6}{x_1} - 20 \leq 0 \tag{42}$$

where $P$= 50,000 N, $E$ = 2×107 N/cm² and the initial design space are: $1 \leq x_i \leq 5$ (i=1,2,..,5), and $30 \leq x_j \leq 65$ (j=6,7,..,10).

BA has achieved a solution that satisfies all the constraints and it reaches the best solution, possibly the unique global optimum. BA outperforms the previous other methods in terms of the minimum objective function value. Table 10 presents the results obtained by BA. We can see that the proposed method requires 25 bats and 1,000 iterations to reach the optimum.

**Table 10** Best solution results for the stepped cantilever beam examples using BA.

This nonlinear constrained problem has been solved by other researchers shown in Table 11. As it is seen, BA significantly outperforms other studies.

**Table 11** Statistical results of the stepped cantilever beam example using different methods.

### 4.7. Case 7: Heat Exchanger Design

As another case study, we now try to solve the heat exchanger design task, which is a difficult benchmark minimization problem since all the constraints are binding. It involves eight design variables and six inequality constraints (three linear and three non-linear). The problem is expressed as follows:

Minimize: $f(X) = x_1 + x_2 + x_3$ (43)

Subject to:

$$g_1 = 0.0025(x_4 + x_6) - 1 \leq 0 \quad (44)$$

$$g_2 = 0.0025(x_5 + x_7 - x_4) - 1 \leq 0 \quad (45)$$

$$g_3 = 0.01(x_8 - x_5) - 1 \leq 0 \quad (46)$$

$$g_4 = 833.33252 x_4 + 100 x_1 - x_1 x_6 - 83333.333 \leq 0 \quad (47)$$

$$g_5 = 1250 x_5 + x_2 x_4 - x_2 x_7 - 125 x_4 \leq 0 \quad (48)$$

$$g_6 = x_3 x_5 - 2500 x_5 - x_3 x_8 + 125 \times 10^4 \leq 0 \quad (49)$$

Table 12 shows the best solution for the heat exchanger design obtained by BA as well as the best solutions obtained previously by other methods. The solution shown for BA is the best generated using 25 Bats. The solution generated by BA (with $X^*$ = [579.30675, 1359.97076, 5109.97052, 182.01770, 295.60118, 217.98230, 286.41653, 395.60118]) is better than the best solutions reported in the literature. As shown in Table 12, the standard deviation and the number of evaluations using BA are also much less than those obtained by the other methods. This solution is feasible and the constraint values are $G^*$ = [0.0000000, 0.0000000, 0.0000000, -0.0071449, -0.0061782, -0.0020000].

**Table 12** Statistical results of the heat exchanger design example by different model.

### 4.8. Case 8: Car Side Impact Design

Design of car side impact is used as a benchmark problem of the proposed BA. On the foundation of European Enhanced Vehicle-Safety Committee (EEVC) procedures, a car is exposed to a side-impact (Youn et al. 2004). Here we want to minimize the weight using nine influence parameters including, thicknesses of B-Pillar inner, B-Pillar reinforcement, floor side inner, cross members, door beam, door beltline reinforcement and roof rail ($x_1$-$x_7$), materials of B-Pillar inner and floor side inner ($x_8$ and $x_9$) and barrier height and hitting position ($x_{10}$ and $x_{11}$). The car side problem is formulated as follow:

Minimize $f(\mathbf{x})$ = Weight;                                                                                                 (50)

Subject to

$g_1(\mathbf{x}) = F_a$ (load in abdomen) ≤ 1 kN;                                                                           (51)

$g_2(\mathbf{x}) = V \times Cu$ (dummy upper chest) ≤ 0.32 m/s;                                                    (52)

$g_3(\mathbf{x}) = V \times Cm$ (dummy middle chest) ≤ 0.32 m/s;                                                 (53)

$g_4(\mathbf{x}) = V \times Cl$ (dummy lower chest) ≤ 0.32 m/s;                                                    (54)

$g_5(\mathbf{x}) = \Delta_{ur}$ (upper rib deflection) ≤ 32 mm;                                                          (55)

$g_6(\mathbf{x}) = \Delta_{mr}$ (middle rib deflection) ≤ 32 mm;                                                         (56)

$g_7(\mathbf{x}) = \Delta_{lr}$ (lower rib deflection) ≤ 32 mm;                                                           (57)

$g_8(\mathbf{x}) = F_p$ (Pubic force) ≤ 4 kN;                                                                                   (58)

$g_9(\mathbf{x}) = V_{MBP}$ (Velocity of V-Pillar at middle point) ≤ 9:9 mm/ms;                          (59)

$g_{10}(\mathbf{x}) = V_{FD}$ (Velocity of front door at V-Pillar) ≤ 15:7 mm/ms;                          (60)

with simple bounds

0.5 ≤ $x_1, x_3, x_4$ ≤ 1.5; 0.45 ≤ $x_2$ ≤ 1.35; 0.875 ≤ $x_5$ ≤ 2.625; 0.4 ≤ $x_6, x_7$ ≤ 1.2; $x_8, x_9 \in$ {0.192, 0.345}; 0.5 ≤ $x_{10}, x_{11}$ ≤ 1.5;

For solving this problem, we ran BA with 20 bats and 1000 iterations. Because this case study has not been solved previously in the literature, we also solved this problem using PSO, DE and GA methods so as to benchmark and compare with the BA method. Table 13 shows the statistical results for the car side impact design problem using the proposed BA method and other well-known methods after 20,000 searches. As it can be seen from Table 13, in comparison with other heuristic algorithms, the proposed algorithm is better than GA and it seems that the BA method performances similar to the PSO and DE.

**Table 13** Statistical results of the car side design example by different methods

## 5. Discussions and Conclusions

We have presented a new bat algorithm for solving engineering optimization problems. BA has been validated using several benchmark engineering design problems, and it is found from our simulations that BA is very efficient. The extensive comparison study, carried out over seven different nonlinear constrained design tasks, reveals that BA performs superior to many different existing algorithms used to solve these seven benchmark problems. It is potentially more powerful than other methods such as GA and PSO as well as harmony search. The primary reason is that BA uses a good combination of major advantages of these algorithms in some way. Moreover, PSO and harmony search are the special cases of BA under appropriate simplifications. More specifically, if we fix the loudness as $A_i=0$ and pulse emission rate as $r_i=1$, BA reduces to the standard particle swarm optimization. On the other hand, if set $A_i=r_i=0.7$ to 0.9, BA essentially becomes a harmony search as frequency change is equivalent to the pitch adjustment in harmony search.

Sensitivity studies can be an important issue for the further research topics, as the fine adjustment of the parameters α and γ can affect the convergence rate of the bat algorithm. This is true for almost all metaheuristic algorithms. In fact, parameter α plays a similar role as the cooling schedule in the simulated annealing. Though the implementation is more complicated than many other metaheuristic algorithms; however, the detailed study of seven engineering design tasks indicates that BA actually uses a balanced combination of the advantages of existing successful algorithms with innovative feature based on the echolocation behavior of microbats. New solutions are generated by adjusting frequencies, loudness and pulse emission rates, while the proposed solution is accepted or not depends on the quality of the solutions controlled or characterized by loudness and pulse rate which are in turn related to the closeness or the fitness of the locations/solution to the global optimal solution.

Theoretically speaking, if we simplify the system with enough approximations, it is possible to analyze the behaviour of the bat algorithm using analysis in the framework of dynamical systems. In addition, more extensive comparison studies with a more wide range of existing algorithms using much tough test functions in higher dimensions will pose more challenges to the algorithms, and thus such comparisons will potentially reveal the virtues and weakness of all the algorithms of interest. Furthermore, a natural extension is to formulate a discrete version of bat algorithm so that it can directly solve combinatorial optimization problems such as the travelling salesman problem. On the other hand, for dynamical optimization problems and computational geometry, a further natural extension to the current bat algorithm would be to use the directional echolocation and Doppler effect, which may lead to even more interesting variants and new algorithms. These further extensions will help us to design more efficient, often hybrid, algorithms to solve a wider class of even tougher optimization problems.